\documentclass[oneside,english]{amsart}
\usepackage[latin1]{inputenc}
\usepackage{amssymb}

\makeatletter

\providecommand{\LyX}{L\kern-.1667em\lower.25em\hbox{Y}\kern-.125emX\@}

 \theoremstyle{plain}
 \newtheorem{thm}{Theorem}[section]
 \numberwithin{equation}{section} 
 \numberwithin{figure}{section} 
 \theoremstyle{remark}
 \newtheorem*{acknowledgement*}{Acknowledgement}
 \theoremstyle{plain}
 \newtheorem{conjecture}[thm]{Conjecture} 
 \theoremstyle{definition}
 \newtheorem{defn}[thm]{Definition}
 \theoremstyle{plain}
 \newtheorem{lem}[thm]{Lemma} 
 \theoremstyle{plain}
 \newtheorem{prop}[thm]{Proposition} 
 \theoremstyle{remark}
 \newtheorem{rem}[thm]{Remark}

\usepackage{babel}
\makeatother

\newcommand{\newcom}{\newcommand}       %
\newcom{\rnewcom}{\renewcommand}
\newcom{\ens}{\ensuremath}
\newcom{\Z}{\mathbb{Z}}                 
\newcom{\N}{\mathbb{N}}                 
\newcom{\Zp}{\Z_p}                      %
\newcom{\pnZ}{p^n\Z}                    %
\newcom{\SLmZp}{SL_m(\Zp)}              %
\newcom{\slm}{\mathfrak{sl}_{_m}}       %
\newcom{\sltwo}{\mathfrak{sl}_{_2}}     %
\newcom{\slmZp}{\slm(\Zp)}              %
\newcom{\gam}{\Gamma_}                  
\newcom{\alf}{\alpha}
\newcom{\bet}{\beta}
\newcom{\stam}[1]{}         
\newcom{\TBD}[1]{}         
\newcom{\middle}[1]{\[ #1 \]}         
\newcom{\rr}{\right}
\rnewcom{\ll}{\left}

\begin{document}

\title{Uniform poly-log diameter bounds for some families of finite groups}

\author{oren dinai}

\begin{abstract}
Fix a prime $p$ and an integer $m$ with $p
> m \geq 2$. Define the family of finite groups
\[ G_{n}:=SL_{m}\left(\mathbb{Z}/p^{n}\mathbb{Z}\right)\]
for $n=1,2,\ldots $. We will prove that there exist two positive
constants $C$ and $d$ such that for any $n$ and any generating set
$S\subseteq G_{n}$,

\[
diam(G_{n},S)\leq C\cdot log^{d}(\left|G_{n}\right|)\] when
$diam\left(G,S\right)$ is the diameter of the finite group $G$
with respect to the set of generators $S$. It is defined as the
maximum over $g\in G$ of the length of the shortest word in $S\cup
S^{-1}$ representing $g$.

This result shows that these families of finite groups have a
poly-logarithmic bound on the diameter with respect to \emph{any}
set of generators. The proof of this result also provides a
efficient algorithm for finding such a poly-logarithmic
representation of any element. In addition it shows that the power
$d$ in the $log$ bound can be arbitrary close to 3 for $m=2$ and
arbitrary close to 4 for $m>2$.

\end{abstract}
\maketitle

\section{Introduction}

The diameter of a finite group $G$ with respect to a set of
generators $S$ is defined to be the diameter of the corresponding
Cayley graph, i.e.,  the minimal number $k$ for which any element
in $G$ can be written as a product of at most $k$ elements in
$S\cup S^{-1}$. We denote this number $diam(G,S)$. We will be
interested in minimizing the diameter of a group with respect to
\emph{any} set of generators. For this we define \[
diam_{worst}(G):=\max _{S\subset G}\left\{ diam(G,S):\textrm{
}S\textrm{ generates }G\right\} \]

While quite a lot is known about the {}``best'' generators, i.e. a
small number of generators which produce a relatively small
diameter (see \cite{BHKLS_diam2}), very little is known about the
worst case.

A well known conjecture of Babai (see \cite{BS1_diam_conj}) asserts:

\begin{conjecture}
\label{conj: Babai}(Babai) There exist constants $d,C$ such that
for any finite non-abelian simple group $G$ \[ diam_{worst}(G)\leq
C\cdot log^{d}(\left|G\right|)\] This bound may even be true for
$d=2$, but not for smaller $d$, as the groups $A_{n}$ demonstrate.
\end{conjecture}
But as of now, there is no family of finite simple groups for
which the Babai's conjecture holds (see \cite{BS1_diam_conj,BS2}
for the best known results).

The goal of this paper is to present for the first time, as far as
we know, a family of finite groups with a poly-logarithmic bound
for the worst-diameter, and to give also an algorithm for
calculating such a poly-logarithmic representation. Our groups are
not simple, though.

\begin{thm}
\label{thm}Fix a prime $p$ such that $p>m\geq 2$ and define
$G_{n}:=SL_{m}\left(\mathbb{Z}/p^{n}\mathbb{Z}\right)$, then for
every real number $d>4$ the following holds:\[
diam_{worst}(G_{n})=O\left(\log ^{d}(\left|G_{n}\right|)\right)\]
Furthermore we will show that if $m>2$, $p$ can be chosen equal to
$m$ and if $m=2$, $d$ can be arbitrary close to 3.

\end{thm}
Our method of proof is a slight improvement of the work of Gamburd
and Shahshahani \cite{GS_gamburd}. Their work was influenced by the
Solovay-Kitaev Lemma (see \cite{NC_solovay_kitaev}).

\section{Preliminaries}

We first restrict ourselves to the case of $m=2$, and then
consider the required modifications for proving the more general
case. From now on we assume that $p$ is an odd prime and $S$ is a
generating set for
$G_{n}:=SL_{m}\left(\mathbb{Z}/p^{n}\mathbb{Z}\right)$. Both are
\emph{arbitrary} but \emph{fixed}. From now $\log{x}$ stands for
$\log_2{x}$ and $\Z_p$ stands for the p-adic integers. First we
begin with some definitions:

\begin{defn} \label{def1}
For any integer $n\geq 0$, define\[ \Gamma _{n}:=\left\{ \gamma
\in SL_{2}\left(\mathbb{Z}\right):\textrm{ }\gamma \equiv
\left(\begin{array}{cc}
 1 & 0\\
 0 & 1\end{array}
\right)(\textrm{mod }p^{n})\textrm{ }\right\} \]

\end{defn}
Equivalently, $\Gamma _{n}=\ker
\left(SL_{2}\left(\mathbb{Z}\right)\stackrel{\pi
_{n}}{\longrightarrow
}SL_{2}\left(\mathbb{Z}/p^{n}\mathbb{Z}\right)\right)$ where $\pi
_{n}$ is the natural projection. Note that with the above
definitions we get that $\Gamma
_{0}=SL_{2}\left(\mathbb{Z}\right)$ and $\Gamma _{n+1}\subset
\Gamma _{n}$ for any $n$. Since $\pi _{n}$ is \emph{onto} we have
\[ G_{n}\cong \Gamma _{0}/\Gamma _{n}\] By abuse of notation,
instead of doing calculations in $G_{n}$ we will do them in
$\Gamma _{0}\textrm{ mod }\Gamma _{n}$, so we will treat the
elements in $G_{n}$ as being in $\Gamma _{0}$.

\begin{defn} \label{def2}
\emph{\label{def: net}}For two subsets $X\subset Y\subset
SL_{2}\left(\mathbb{Z}\right)$ denote $X\equiv Y\left(\textrm{mod
}p^{k}\right)$ if $\pi _{k}\left(X\right)=\pi _{k}\left(Y\right)$.
Denote also\begin{equation} Y\stackrel{l}{\rightsquigarrow
}X\label{eq:Y l X}\end{equation} if every element in $Y$ can be
moved into $X$ by a multiplication of at most of $l$ elements in
$S\cup S^{-1}$. Explicitly, $\forall  y\in Y$, $\exists
s_1,s_2,\ldots,s_k\in S\cup S^{-1}$, for some $k\leq l$, such that
$y\cdot s_1\cdot s_2\cdot\ldots\cdot s_k\in X$.
\end{defn}
We need to distinguish between the group commutator and the Lie
bracket operations. For elements $g,h$ in a group $G$ we denote
$\left\{ g,h\right\} :=g^{-1}h^{-1}gh$ for the group commutator.
For elements $A,B$ in the Lie algebra
$\mathfrak{sl}_{_{m}}\left(R\right)$ for some commutative ring $R$
we write $\left[A,B\right]:=AB-BA$ the Lie bracket, when
$\mathfrak{sl}_{_{m}}\left(R\right)$ is the set of
 matrices in $M_{m}(R)$ with Trace $=0$.

\section{statement of the main results}

The first Lemma \ref{lemAlex} concerns with some generation
properties in the Lie algebras
$\mathfrak{sl}_{_{m}}\left(R\right)$. The next three statements
\ref{lem},\ref{prop} and \ref{thm:1'} are restricted to the case
$m=2$ while Theorem \ref{thm:generalCase} is a generalization of
these statements to $m\geq 2$. Lemma \ref{lem} is due to Michael
Larsen. We will see later that this Lemma already has almost all
the key ideas for proving Theorem \ref{thm}. This Lemma will give
us a reduction from the groups $SL_m(\Z_p)$ to the algebras
$\mathfrak{sl}_{_{m}}(\Z_p)$.

The following Lemma is the main ingredient for Theorem
\ref{thm:generalCase} which is a generalization of \ref{thm:1'}.

\begin{lem}
\label{lemAlex}For any prime $p$ and integer $m$ with $p>m\geq 2$,
any element in $\mathfrak{sl}_{_{m}}\left(\mathbb{Z}_{p}\right)$
is a sum of two Lie brackets. This is true also for $p\geq m> 2$.
Furthermore, in the case $p>m=2$, any element in
$\mathfrak{sl}_{_{2}}\left(\mathbb{Z}_{p}\right)$ can be expressed
as one Lie bracket.
\end{lem}

The following Lemma \ref{lem} and Proposition \ref{prop} are
restricted to the case $m=2$. Remember that we have defined
$\Gamma_n=\left\{ \gamma \in
SL_{2}\left(\mathbb{Z}\right):\textrm{ }\gamma \equiv
\left(\begin{array}{cc}
 1 & 0\\
 0 & 1\end{array}
\right)(\textrm{mod }p^{n})\textrm{ }\right\}$.

\begin{lem} \label{lem}
For any integers $i,j\geq 0$ and for any $k\leq \min
\left\{ i,j\right\} $ the group commutator map
\begin{eqnarray}
\Gamma _{i}/\Gamma _{i+k}\times \Gamma _{j}/\Gamma _{j+k} & \stackrel{\left\{ \cdot ,\cdot \right\} }{\longrightarrow } & \Gamma _{i+j}/\Gamma _{i+j+k}\label{eq:commutator map}\\
\left(\overline{\alpha },\overline{\beta }\right) & \mapsto  & \overline{\left\{ \alpha ,\beta \right\} }\nonumber
\end{eqnarray}
is surjective.
\end{lem}
Lemma \ref{lem} will imply the following iteration step needed for
 Theorem \ref{thm:1'}:

\begin{prop} \label{prop}
For any $d>2$ there exists $C$ such that for any $n\geq 1$
\[\Gamma _{n}\stackrel{Cn^{d}}{\rightsquigarrow }\Gamma _{n+1} \]

\end{prop}
With this proposition we will prove Theorem \ref{thm:1'}, which is
equivalent to Theorem \ref{thm} for the case $m=2$.

\begin{thm}
\label{thm:1'}Fix $p>2$ and set
$G_{n}:=SL_{2}\left(\mathbb{Z}/p^{n}\mathbb{Z}\right)$ then for
any $d>3$ there exists a real constant $C$ such that for any set
of generators $S\subseteq G_{n}$, any element in $G_{n}$ can be
written as a product of at most $Cn^{d}$ elements in $S\cup
S^{-1}$.
\end{thm}

The next Theorem is equivalent to Theorem \ref{thm} for the case
$m>2$ (see Remark \ref{remark}).

\begin{thm}
\label{thm:generalCase}Fix a prime $p$ and an integer $m$ with
$p\geq m>2$. Define the family of finite groups
$G_{n}:=SL_{m}\left(\mathbb{Z}/p^{n}\mathbb{Z}\right)$  for
$n=1,2,\ldots $. For any real $d$ with $d>4$ there exists a real
constant $C$ such that for any set of generators $S\subseteq
G_{n}$, any element in $G_{n}$ can be written as a product of at
most $Cn^{d}$ elements in $S\cup S^{-1}$.
\end{thm}

\section{proofs}

\textbf{\emph{Proof of Lemma \ref{lemAlex}:}} Let us denote
$diag(\lambda_1,\ldots,\lambda_m)$ to be the diagonal matrix with
these values on its diagonal. For a matrix
$A\in\mathfrak{sl}_{_{m}}\left(\mathbb{Z}_{p}\right)$ denote by
$diagA$ the diagonal matrix with the same diagonal of $A$. Now
choose $D=diag(\lambda_{1},\ldots ,\lambda_m)$ such that
$\sum_{i=1}^{m}{\lambda_i}=0$ and $(\lambda_i-\lambda_j)$ is a
unit in $ \mathbb{Z}_p$ for any $i \neq j$. Take for example
$\{\lambda_i\}$ to be $\pm 1,\ldots,\pm k$ if $m=2k$ or $0,\pm
1,\ldots,\pm k$ if $m=2k+1$. Since for any $i\neq j$,
$[D,E_{i,j}]=(\lambda_i-\lambda_j)E_{i,j}$ all we need to show is
that given $A\in\mathfrak{sl}_{_{m}}\left(\mathbb{Z}_{p}\right)$
we can find two matrices $B',B''$ s.t. $diag[B',B'']=diagA$. For
if we write $A-[B',B'']=\sum_{i\neq j}{a_{i,j}E_{i,j}}$ we get
$A=[B',B'']+\sum_{i\neq j}{a_{i,j}E_{i,j}}=[B',B'']+[D,\sum_{i
\neq j}{\frac{a_{i,j}}{(\lambda_i-\lambda_j)}}E_{i,j}]$.

Let us denote by $B^g$ the representation of the matrix $B$ in the
basis $e_1,e_1+e_2,e_1+e_2+e_3,\ldots,e_1+\ldots +e_m$ where
$\{e_1,\ldots,e_m\}$ is the standard basis. For any $i<m$ we get
$diag[D^g,E_{i,i+1}^g]=diag[D,E_{i,i+1}]^g=diag(\lambda_i-\lambda_{i+1})E_{i,i+1}^g=(\lambda_i-\lambda_{i+1})(E_{i+1,i+1}-E_{i,i})$.
So if we write $diagA=\sum_{i=1}^{m-1}a_i(E_{i+1,i+1}-E_{i,i})$ we
get that
$diagA=\sum_{i=1}^{m-1}a_i(E_{i+1,i+1}-E_{i,i})=diag[D^g,\sum_{i=1}^{m-1}\frac{a_i}{(\lambda_i-\lambda_{i+1})}E_{i,i+1}^g]$,
so we are done. Note that when $m>2$, $p$ is odd we can take $m$
to be equal $p=2k+1$ and follow the same arguments.

The following improvement to the case $\sltwo(\Z_p)$ is due to
Larsen. In $\sltwo(\Z_p)$ we have the following identity: for
every $C,\textrm{ }D\in
\mathfrak{sl}_{_{2}}\left(\mathbb{Z}_{p}\right)$
\begin{equation}
\left[\left[C,\, D\right],\,
C\right]=2Tr\left(CD\right)C-2Tr\left(C^{2}\right)D\label{eq:1
basic lie identity}\end{equation} (this identity can be checked by
expressing the matrices $C$ and $D$ explicitly via their entries).
From identity \ref{eq:1 basic lie identity} we get the following
identity for every $A,\textrm{ }B\in
\mathfrak{sl}_{_{2}}\left(\mathbb{Z}_{p}\right)$\begin{equation}
\left[\left[\left[A,\, B\right],\, A\right],\, \left[A,\,
B\right]\right]=-2Tr\left(\left[A,\,
B\right]^{2}\right)A\label{eq:3 bracket lie
identity}\end{equation} by setting $C=\left[A,B\right]\textrm{,
}D=A$ and observing that the first term of the right-hand side of
\ref{eq:1 basic lie identity} vanished since
$Tr\left(CD\right)=Tr\left(\left[A,B\right]A\right)=Tr\left(\left[A,BA\right]\right)=0$.

First we want to use identity \ref{eq:3 bracket lie identity}
to show that any element in
$\mathfrak{sl}_{_{2}}\left(\mathbb{Z}_{p}\right)\setminus p\cdot
\mathfrak{sl}_{_{2}}\left(\mathbb{Z}_{p}\right)$ is a bracket. If
$A=\left(\begin{array}{cc}
 u & v\\
 w & -u\end{array}
\right)$ is not in $p\cdot
\mathfrak{sl}_{_{2}}\left(\mathbb{Z}_{p}\right)$ it has at least
one entry which is unit in $\mathbb{Z}_{p}$. By a straightforward
calculation we get that for $B=\left(\begin{array}{cc}
 0 & 1\\
 0 & 0\end{array}
\right),\left(\begin{array}{cc}
 0 & 0\\
 1 & 0\end{array}
\right)\textrm{ or }\left(\begin{array}{cc}
 0 & 1\\
 1 & 0\end{array}
\right)$, $Tr\left(\left[A,\, B\right]^{2}\right)$ equals
$2w^{2},2v^{2}\textrm{ or }-8u^{2}$ respectively. Therefore
$-2\beta Tr\left(\left[A,\, B\right]^{2}\right)$ equals $1$ for
some $B\in \mathfrak{sl}_{_{2}}\left(\mathbb{Z}_{p}\right)\textrm{
and some }\beta \in \mathbb{Z}_{p}$. Note that we used here the
fact that $p$ is odd and so $2$ is unit in $\mathbb{Z}_{p}$. So we
can find $A_{1},A_{2}\in
\mathfrak{sl}_{_{2}}\left(\mathbb{Z}_{p}\right)$ when $A_{1}=\beta
\left[\left[A,\, B\right],\, A\right]\textrm{, }A_{2}=\left[A,\,
B\right]$ s.t. $\left[A_{1},A_{2}\right]=\left[\beta
\left[\left[A,\, B\right],\, A\right],\, \left[A,\,
B\right]\right]=A$ as we wanted.

Now we show that any element $A$ in $\sltwo(\Z_p)$ can be
expressed as one Lie bracket. For $A=0$ the statement is clear so
take $A\neq 0$ in $\sltwo(\Z_p)$, then there exists $k$ such that
$A'=p^{-k}A$ and $A'\in\sltwo(\Z_p)\setminus p\cdot \sltwo(\Z_p)$.
By the previous paragraph there exist $B',B''$ with $A'=[B',B'']$
and so $A=[p^kB',B'']$ and we are done.
\begin{flushright} $\square$ \end{flushright}

\textbf{\emph{Proof of Lemma \ref{lem}:}} First we observe that
the commutator map \ref{eq:commutator map} is well defined since
for any $\alpha \in \Gamma _{i}\, ,\, \beta \in \Gamma _{j}\, ,\,
\textrm{ }\alpha '\in \Gamma _{i+k}\, ,\, \beta '\in \Gamma
_{j+k}$ there exists $4$ matrices $A,B,A',B'\in
SL_{2}\left(\mathbb{Z}\right)$ s.t. $\alpha =I+p^{i}A\, ,\, \alpha
'=I+p^{i+k}A'\, ,\, \beta =I+p^{j}B\, ,\, \beta '=I+p^{j+k}B'$ and
so we get \[ \left\{ \alpha ,\beta \right\} \equiv \left\{ \alpha
\alpha ',\beta \beta '\right\} \equiv
I+p^{i+j}\left[A,B\right]\left(\textrm{mod }p^{i+j+k}\right)\]
Secondly we observe that we can work $p$-adicly which means that
we can do all the calculations over $\mathbb{Z}_{p}$-the ring of
$p$-adic integers instead over $\mathbb{Z}$. Indeed if we denote
$\overline{\Gamma _{m}}:=\ker
\left(SL_{2}\left(\mathbb{Z}_{p}\right)\stackrel{\overline{\pi
_{m}}}{\longrightarrow
}SL_{2}\left(\mathbb{Z}_{p}/p^{m}\mathbb{Z}_{p}\right)\right)$
then we get for any $m,k\geq 0$ \[ \Gamma _{m}/\Gamma _{m+k}\cong
\overline{\Gamma _{m}}/\overline{\Gamma _{m+k}}\] Instead of doing
the group commutator we want to reduce our problem to Lie algebras
and their bracket product which is easier to handle. We have the
following bijections for any $m\geq 1$:\begin{eqnarray}
\overline{\Gamma _{m}} & \stackrel{\log \left(1+x\right)}{\longrightarrow } & p^{m}\mathfrak{sl}_{_{2}}\left(\mathbb{Z}_{p}\right)\label{def: log(1+x)}\\
I+p^{m}A & \mapsto  & \log
\left(I+p^{m}A\right)=p^{m}A-\frac{p^{2m}}{2}A^{2}+\frac{p^{3m}}{3}A^{3}-\ldots
\label{map: log(1+x)}
\end{eqnarray}
i.e. $\log \left(1+x\right):=\sum _{n=1}^{\infty
}\left(-1\right)^{n+1}\frac{x^{n}}{n}$. If we denote by
$\overline{L}:=\mathfrak{sl}_{_{2}}\left(\mathbb{Z}_{p}\right)$
the Lie algebra over $\mathbb{Z}_{p}$ then we get the following
commutative diagram\begin{equation}
\begin{array}{ccc}
 \left(\overline{\Gamma _{i}}/\overline{\Gamma _{i+k}}\right)\times \left(\overline{\Gamma _{j}}/\overline{\Gamma _{j+k}}\right) & \stackrel{\left\{ \cdot ,\cdot \right\} }{\longrightarrow } & \overline{\Gamma _{i+j}}/\overline{\Gamma _{i+j+k}}\\
 \varphi _{1}\downarrow  &  & \downarrow \varphi _{2}\\
 \left(p^{i}\overline{L}/p^{i+k}\overline{L}\right)\times \left(p^{j}\overline{L}/p^{j+k}\overline{L}\right) & \stackrel{\left[\cdot ,\cdot \right]}{\longrightarrow } & p^{i+j}\overline{L}/p^{i+j+k}\overline{L}\end{array}
\label{cummotative diagram}\end{equation} When $\varphi
_{1},\varphi _{2}$ are the bijections $\varphi
_{1}\left(x_{1},x_{2}\right):=\log \left(1+x_{1}\right)\times \log
\left(1+x_{2}\right)\textrm{ and }\varphi _{2}\left(x\right):=\log
\left(1+x\right)$ as in \ref{def: log(1+x)}. Note that the
summation in \ref{map: log(1+x)} indeed converge since we are over
$\mathbb{Z}_{p}$ and its general element converge to zero (for
more details see \cite{dinai}).

\label{paragraph} In order to show that the commutator map
$\left\{ \cdot ,\cdot \right\} $ in \ref{cummotative diagram} is
onto it is enough to show that the bracket map $\left[\cdot ,\cdot
\right]$ in \ref{cummotative diagram} is onto. So it suffices to
show that for every $k\geq 1$, every element in the Lie algebra
$\overline{L}/p^{k}\overline{L}=\mathfrak{sl}_{_{2}}\left(\mathbb{Z}_{p}\right)/p^{k}\mathfrak{sl}_{_{2}}\left(\mathbb{Z}_{p}\right)$
is a bracket of two elements in
$\mathfrak{sl}_{_{2}}\left(\mathbb{Z}_{p}\right)/p^{k}\mathfrak{sl}_{_{2}}\left(\mathbb{Z}_{p}\right)$.
\stam{ Now it is enough to work in
$\mathfrak{sl}_{_{2}}\left(\mathbb{Z}_{p}\right)$ and show that
every element $A'$ in
$\mathfrak{sl}_{_{2}}\left(\mathbb{Z}_{p}\right)\setminus
p^{k}\mathfrak{sl}_{_{2}}\left(\mathbb{Z}_{p}\right)$ is a
bracket. Because the Lie bracket is bi-linear and since there is
$0\leq l\leq k$ and $A\in
\mathfrak{sl}_{_{2}}\left(\mathbb{Z}_{p}\right)\setminus p\cdot
\mathfrak{sl}_{_{2}}\left(\mathbb{Z}_{p}\right)$ s.t $A'=p^{l}A$
it is sufficient to show that every $A$ in
$\mathfrak{sl}_{_{2}}\left(\mathbb{Z}_{p}\right)\setminus p\cdot
\mathfrak{sl}_{_{2}}\left(\mathbb{Z}_{p}\right)$ is a bracket.} By
Lemma \ref{lemAlex} we are done.
\begin{flushright} $\square$ \end{flushright}

\stam{Now by using these Lemmas we get the iteration proposition
we wanted:}

\textbf{\emph{Proof of Proposition \ref{prop}:} }We need to prove
that any element in $\Gamma _{n}/\Gamma _{n+1}$ can be written as
a product in at most $Cn^{d}$ elements in $S\cup S^{-1}$. Let us
denote the minimal length of $\gamma \in \Gamma _{n}/\Gamma
_{n+1}$ by $l\left(\gamma \right)$. We prove that $l\left(\gamma
\right)\leq Cn^{d}$ by induction on $n$. For any $d>2$ we can
choose $N_{0}$ s.t for any $n>N_{0}$,
$4\left(\frac{n+1}{2n}\right)^{d}<1$. Choose a constant $C$ big
enough s.t. $l\left(\gamma \right)\leq Cn^{d}$ for any $\gamma \in
\Gamma _{n}/\Gamma _{n+1}$ and any $n\leq N_{0}$. Now let
$n>N_{0}$ and let $\gamma \in \Gamma _{n}/\Gamma _{n+1}$. There
are always $k,m\leq \frac{n+1}{2}$ with $k+m=n$. Hence by Lemma
\ref{lem} there exists $\gamma _{1}\in \Gamma _{m}/\Gamma _{m+1}$
and $\gamma _{2}\in \Gamma _{k}/\Gamma _{k+1}$ with $\gamma
=\left\{ \gamma _{1},\gamma _{2}\right\} $ and so $l\left(\gamma
\right)\leq 2\left(l\left(\gamma _{1}\right)+l\left(\gamma
_{2}\right)\right)$ and by the induction hypothesis we
get\begin{eqnarray*} l\left(\gamma \right) & \leq  &
2(Ck^{d}+Cm^{d})\leq
4C\left(\frac{n+1}{2}\right)^{d}=4\left(\frac{n+1}{2n}\right)^{d}Cn^{d}<Cn^{d}
\end{eqnarray*}
as claimed. \TBD{ The proof for the case $p>m=2$ was given in the
second half of the proof for \ref{lem}.}
\begin{flushright} $\square$ \end{flushright}

Now it remains to combine the proceeding and get Theorem
\ref{thm:1'} but before that we remark on the equivalence of
Theorems \ref{thm},\ref{thm:1'} and \ref{thm:generalCase}.

\begin{rem} \label{remark}
For the equivalence between the Theorems \ref{thm}, \ref{thm:1'}
and \ref{thm:generalCase} we note that
$\left|G_{n}\right|=\left|SL_{m}\left(\mathbb{Z}/p^{n}\mathbb{Z}\right)\right|=
p^{\theta(n  m^{2})}$. Therefore, $\log \left|G_{n}\right|\sim
\log(p)m^{2}n$ and so the equivalence is clear.
\end{rem}
\textbf{\emph{Proof of Theorem \ref{thm:1'}:}} By applying $(n-1)$
times proposition \ref{prop} and then combining those steps
together we get that for any $d>2$ there exists $C$ such that\[
\Gamma _{1}\stackrel{C(1^{d}+2^{d}+\ldots
+\left(n-1\right)^{d})}{\rightsquigarrow }\Gamma _{n}\textrm{ }\]
If we choose $l_{0}$ such that $\Gamma
_{0}\stackrel{l_{0}}{\rightsquigarrow }\Gamma _{1}$ and we can
assume that $l_{0}\leq C$ then we get $\Gamma
_{0}\stackrel{Cn^{d+1}}{\rightsquigarrow }\Gamma _{n}$. Therefore
we got the result we wanted, for any $d>3$ there exist $C$ such
that\[ diam(SL_{2}\left(\mathbb{Z}/p^{n}\mathbb{Z}\right),S)\leq
Cn^{d}\]
\begin{flushright} $\square$ \end{flushright}

\textbf{\emph{Proof of Theorem \ref{thm:generalCase}:}}  Now we
make the required modification to the definitions and statements
\ref{def1}, \ref{def2}, \ref{lem} and \ref{prop}. Let's replace in
these definitions all the occurrences of $SL_2(),\sltwo()$ by
$SL_m(),\slm()$ respectively.

Let's modify Lemma \ref{lem} to the following: every element in
$\Gamma _{i+j}/\Gamma _{i+j+k}$ can be expressed as a product of
two commutators $\{\alf,\bet\}\{\alf',\bet'\}$ when $\alpha ,
\alpha ' \in \gam i/\gam {i+k}$ and $\beta ,\beta ' \in \gam
j/\gam {j+k}$. The proof of this follow the same lines of the
proof of \ref{lem} but instead of representing every element in
$\sltwo$ by one Lie bracket \stam{(the second half of proof
\ref{lem})}we use the previous Lemma \ref{lemAlex} about
representing each element by a sum of two Lie bracket and so we
get the required representation as product of two commutators.

\TBD{(\textbf{Alex:} here i'm not sure why the \textbf{sum} of two
brackets imply on the \textbf{product} of two commutators, i.e.
why with the operators +,$\cdot$ the diagram commutes. From the
point of commutation of $\cdot$ it works since the groups are
abelian). }

Now lets modify Proposition \ref{prop} to the same claim for any
$d>3$. In its proof we see that if we use the previous Lemma about
expressions as a product of $b$ commutators (we proved it for
$b=2$) then we get $l\left(\gamma \right)\leq
2b\left(l\left(\gamma _{1}\right)+l\left(\gamma
_{2}\right)\right)$ and so
\begin{eqnarray*} l(\gamma)\leq
2b(Ck^{d}+Cm^{d})\leq
4bC(\frac{n+1}{2})^{d}=4b(\frac{n+1}{2n})^{d}Cn^{d}<Cn^{d}
\end{eqnarray*} when the last inequality is true if
$d>\log(4b)=2+\log(b)$ and $n$ is big enough. Now for $b=2$ we get
the result we wanted for any $d>3$.

In conclusion if we combine all the previous modifications we can
use them in the proof of Theorem \ref{thm:1'} to get the
generalization we wanted: for any $d>4$ there exist $C$ such that
\[ diam(SL_{m}\left(\mathbb{Z}/p^{n}\mathbb{Z}\right),S)\leq
Cn^{d}\] for \emph{any} generating set $S$ of $SL_m(\Z/\pnZ)$.
\begin{flushright} $\square$ \end{flushright}

\begin{acknowledgement*}
This work is part of the author's M.sc thesis. I wish to thank my
advisor Prof. Alex Lubotzky for sharing his ideas, for his
guidance and assistance in writing this work. I wish to thank also
Prof. Michael Larsen for Lemma \ref{lem} and to the referee for
his careful report and for many of his suggestions and
observations.
\end{acknowledgement*}


\begin{thebibliography}{BHKLS}
\bibitem[BKL]{BKL_diam1}Babai, L., Kantor, W.M., Lubotzky, A.: Small diameter Cayley graphs
for finite simple groups, Europ.J.Combinatorics, 10, (1989), 507-522.
\bibitem[BHKLS]{BHKLS_diam2}Babai, L., Hetyei, G., Kantor, W. M., Lubotzky, A., Seress, A.: On
the diameter of finite groups. In 31st Annual Symposium on Foundations
of Computer Science, volume II, pages 857-865, St. Louis, Missouri,
22-24 October 1990. IEEE.
\bibitem[BS1]{BS1_diam_conj}Babai, L., Seress, A.: On the diameter of Cayley graphs of the symmetric
group, J.Combinatorial Theory-A 49 , (1988), 175-179.
\bibitem[BS2]{BS2}Babai, L., Seress, A.: On the diameter of permutation groups, Europ.
J. Comb. 13, (1992), 231-243.
\bibitem[Di]{dinai}Dinai, O.: Poly-log diameter bounds for some families of finite groups, Master's thesis, Hebrew University, (2004).
\bibitem[GS]{GS_gamburd}$\textrm{Gamburd,A.,Shahshahani}$, M.: Uniform diameter bounds for
some families of Cayley graphs, Internat. Math. Res. Notices, 71,
(2004), 3813-3824.
\bibitem[NC]{NC_solovay_kitaev}M.A.Nielsen, I.L.Cuang, Quantum computation and quantum information,
Cambridge University Press, Cambridge, (2000).\end{thebibliography}
\end{document}